\input amstex
\input amsppt.sty
\magnification=\magstep1
\hsize=30truecc
\vsize=22.2truecm
\baselineskip=16truept
\nologo
\pageno=1
\TagsOnRight
\topmatter
\def\Z{\Bbb Z}
\def\N{\Bbb N}

\def\C{\Bbb C}

\def\M{\Bbb M}
\def\l{\left}
\def\r{\right}
\def\bg{\bigg}
\def\({\bg(}
\def\[{\bg[}
\def\){\bg)}
\def\]{\bg]}
\def\t{\text}
\def\f{\frac}

\def\bi{\binom}

\def\ls{\leqslant}
\def\gs{\geqslant}
\def\al{\alpha}

\def\da{\delta}

\def\bi{\binom}

\def\Proof{\noindent{\it Proof}}
\def\Remark{\medskip\noindent{\it  Remark}}
\def\Ack{\medskip\noindent {\bf Acknowledgment}}

\def\M#1#2#3#4{\left[\matrix #1&#2\\#3&#4\endmatrix\right]_n}
\def\m#1#2#3#4{\left[\matrix #1&#2\\#3&#4\endmatrix\right]_{n-1}}
\def\mm#1#2#3#4{\left[\matrix #1&#2\\#3&#4\endmatrix\right]}
\topmatter \hbox{Acta Arith. 125(2006), no.\,1, 21--39.}
\bigskip
\title {Identities concerning Bernoulli and Euler polynomials}\endtitle
\author {Zhi-Wei Sun and Hao Pan (Nanjing)}\endauthor
\leftheadtext{Z. W. Sun and H. Pan}
\abstract We establish two
general identities for Bernoulli and Euler polynomials, which are
of a new type and have many consequences. The most striking result
in this paper is as follows: If $n$ is a positive integer,
$r+s+t=n$ and $x+y+z=1$, then we have
$$r\M stxy+s\M tryz+t\M rszx=0$$
where
$$\M stxy:=\sum_{k=0}^n(-1)^k\bi sk\bi t{n-k}B_{n-k}(x)B_k(y).$$
It is interesting to compare this with the following property
of determinants:
$$r\left|\matrix s&t\\x&y\endmatrix\right|
+s\left|\matrix t&r\\y&z\endmatrix\right|
+t\left|\matrix r&s\\z&x\endmatrix\right|=0.$$
Our symmetric relation implies the curious identities of Miki and Matiyasevich as well as
some new ones for Bernoulli polynomials such as
$$\sum_{k=0}^n{\bi nk}^2B_k(x)B_{n-k}(x)=2\sum^n\Sb k=0\\k\not=n-1\endSb
\bi nk\bi {n+k-1}kB_k(x)B_{n-k}.$$
\endabstract
\thanks  2000 {\it Mathematics Subject Classification}.
Primary 11B68; Secondary 05A19.
\newline\indent The first author is responsible for communications,
and partially supported by the National Science Fund for
Distinguished Young Scholars (no. 10425103) and a Key Program of
NSF (no. 10331020) in China.
\endthanks
\endtopmatter
\document
\hsize=30truecc
\vsize=22.2truecm
\baselineskip=16truept

\heading 1. Introduction\endheading
Let $\N=\{0,1,2,\ldots\}$ and $\Z^+=\{1,2,3,\ldots\}$.
The well-known Bernoulli numbers $B_n\ (n\in\N)$ are rational numbers defined by
$$B_0=1\ \ \t{and}\ \ \sum_{k=0}^n\bi{n+1}kB_k=0\ \ (n\in\Z^+).$$
Similarly, Euler numbers $E_n\ (n\in\N)$ are integers given by
$$E_0=1\ \ \t{and}\ \ \sum^n\Sb k=0\\2\mid n-k\endSb\bi nkE_k=0\ \ (n\in\Z^+).$$

For $n\in\N$ the Bernoulli polynomial $B_n(x)$ and the Euler polynomial $E_n(x)$ are
as follows:
$$B_n(x)=\sum_{k=0}^n\bi nkB_kx^{n-k}\ \ \t{and}\ \ E_n(x)
=\sum_{k=0}^n\bi nk\f {E_k}{2^k}\l(x-\f12\r)^{n-k}.$$
Clearly $B_n(0)=B_n$ and $E_n(1/2)=E_n/2^n$.
Here are some basic properties of Bernoulli and Euler polynomials
we will need later.
$$\gather B_n(1-x)=(-1)^nB_n(x),\ \Delta (B_n(x))=nx^{n-1};
\\E_n(1-x)=(-1)^nE_n(x),\ \Delta^*(E_n(x))=2x^n.
\endgather$$
(The operators $\Delta$ and $\Delta^*$ are defined by $\Delta(f(x))=f(x+1)-f(x)$
and $\Delta^*(f(x))=f(x+1)+f(x)$.)
 Also, $B_{n+1}'(x)=(n+1)B_n(x)$ and $E_{n+1}'(x)=(n+1)E_n(x)$.

For a sequence $\{a_n\}_{n\in\N}$ of complex numbers, its dual sequence $\{a_n\}_{n\in\N}$ is
given by $a_n^*=\sum_{k=0}^n\bi nk(-1)^ka_k\ (n\in\N)$. It is well known that $a_n^{**}=a_n$.
In 2003 Z. W. Sun [S2] deduced some combinatorial identities in dual sequences.
The sequences $\{(-1)^nB_n\}_{n\in\N}$ and $\{(-1)^nE_n(0)\}_{n\in\N}$
are both self-dual sequences (cf. [S2]),
later we will make use of this fact.

In 1978 H. Miki [Mi] discovered the following curious identity:
$$\sum_{k=2}^{n-2}\f{B_kB_{n-k}}{k(n-k)}
-\sum_{k=2}^{n-2}\bi nk\f{B_kB_{n-k}}{k(n-k)}=2H_n\f{B_n}n$$
for every $n=4,5,\ldots$, where $H_n=1+1/2+\cdots+1/n$.
In 1997 Y. Matiyasevich [Ma] found another
identity of this type:
$$(n+2)\sum_{k=2}^{n-2}B_kB_{n-k}-2\sum_{k=2}^{n-2}\bi{n+2}kB_kB_{n-k}=n(n+1)B_n$$
for any $n=4,5,\ldots$. These two identities
are of a deep nature.
In fact, all known proofs of these identities by others are complicated
(cf. [Mi], [G] and [DS]); for example,
the approach of G. V. Dunne and C. Schubert [DS] was even motivated
by quantum field theory and string theory.

Recently the authors [PS] presented a new method to handle such
identities. Though their approach only involves differences and
derivatives of polynomials, they were able to use the powerful
method to extend Miki's and Matiyasevich's identities to
identities concerning the sums $\sum_{k=0}^nB_k(x)B_{n-k}(y)$ and
$$\sum_{k=1}^{n-1}\f {B_k(x)}k\cdot\f{B_{n-k}(y)}{n-k}
=\f1n\sum_{k=1}^{n-1}\f {B_k(x)}kB_{n-k}(y)
+\f1n\sum_{l=1}^{n-1}\f{B_l(y)}lB_{n-l}(x)$$
(where $n$ is a positive integer). They also handled similar sums related to Euler polynomials.

Let $n$ be any positive integer. As usual, $\bi zn=z(z-1)\cdots(z-n+1)/n!$
(and $\bi z0=1$) even if
$z\not\in\N$. Observe that
$$\sum_{k=0}^nB_k(x)B_{n-k}(y)=\sum_{k=0}^n(-1)^k\bi{-1}kB_k(x)B_{n-k}(y)$$
and
$$\align-\sum_{k=1}^n\f{B_k(x)}kB_{n-k}(y)
&=\sum_{k=1}^n(-1)^{k}\bi{-1}{k-1}\f{B_k(x)}kB_{n-k}(y)
\\&=\lim_{t\to0}\f1t\sum_{k=1}^n(-1)^k\bi tkB_k(x)B_{n-k}(y).
\endalign$$
Inspired by this observation, here we investigate relations among the sums
$$\sum_{k=0}^n(-1)^k\bi sk\bi t{n-k} P_k(x)Q_{n-k}(y)$$
with $P,Q\in\{B,E\}$.

 Our central result is the following theorem.

\proclaim{Theorem 1.1} Let $n\in\Z^+$ and $x+y+z=1$.

{\rm (i)} If $r+s+t=n-1$, then
$$\aligned&\sum_{k=0}^n(-1)^k\bi rk\bi s{n-k}B_k(x)E_{n-k}(z)
\\&-(-1)^n\sum_{k=0}^n(-1)^k\bi rk\bi t{n-k}B_k(y)E_{n-k}(z)
\\=&\f r2\sum_{l=0}^{n-1}(-1)^l\bi sl\bi t{n-1-l}E_l(y)E_{n-1-l}(x).
\endaligned\tag1.1$$

{\rm (ii)}
If $r+s+t=n$, then we have the symmetric relation
$$r\M stxy+s\M tryz+t\M rszx=0\tag1.2$$
where
$$\M stxy:=\sum_{k=0}^n(-1)^k\bi sk\bi t{n-k}B_{n-k}(x)B_k(y).\tag1.3$$
\endproclaim

\Remark\ 1.1. It is interesting to compare (1.2) with the following property
of determinants:
$$0=\left|\matrix r&s&t\\r&s&t\\z&x&y\endmatrix\right|
=r\left|\matrix s&t\\x&y\endmatrix\right| +s\left|\matrix
t&r\\y&z\endmatrix\right| +t\left|\matrix
r&s\\z&x\endmatrix\right|.$$ In view of K. Dilcher's paper [D],
the referee thought that Theorem 1.1 might have a generalization
involving sums of products of $m$ Bernoulli or Euler polynomials.
But we are unable to obtain a compact extension of Theorem 1.1
though we have made a serious attempt.

\proclaim{Corollary 1.1} Let $n\in\Z^+$ and let $\al,x,y$ be parameters. Then
$$\aligned&\f{\al+n+1}2\sum_{k=0}^{n-1}\bi{\al+k}kE_k(x)E_{n-1-k}(y)
\\=&\sum_{k=0}^n\bi{\al+n+1}k\l((-1)^{n-k}B_k(x)-\bi{\al+n-k}{n-k}B_k(y)\r)E_{n-k}(x-y)
\endaligned\tag1.4$$
 and
$$\aligned&(\al+n+2)\sum_{k=0}^n\bi{\al+k}kB_k(x)B_{n-k}(y)
\\=&(\al+1)\sum_{k=0}^n\bi{\al+n+2}k(-1)^{n-k}B_k(x)B_{n-k}(x-y)
\\&+\sum_{k=0}^n\bi{\al+n+2}k\bi{\al+n-k}{n-k}B_k(y)B_{n-k}(x-y).
\endaligned\tag1.5$$
\endproclaim
\Proof. Let $x'=1-x$ and $z'=x-y$. Then $x'+y+z'=1$.
Applying Theorem 1.1(i) with $r=\al+n+1$, $s=-1$ and $t=-\al-1$ we then get (1.4).
(Note that $(-1)^k\bi{-z}k=\bi{z+k-1}k$.)
By Theorem 1.1(ii),
$$\align&(\al+n+2)\M{-1}{-\al-1}{1-x}y
\\=&\M{-\al-1}{\al+n+2}y{x-y}+(\al+1)\M{\al+n+2}{-1}{x-y}{1-x}.
\endalign$$
This is an equivalent version of (1.5). \qed

\Remark\ 1.2. (1.5) in the case $\al=x=y=0$ yields Matiyasevich's identity
since $B_{2l+1}=0$ for $l=1,2,3,\ldots$.

\proclaim{Corollary 1.2} Let $n>l\gs0$ be integers. Then
$$\aligned&\f{n-l+1}2\sum_{k=\da_{l,0}}^{n}\bi nk\bi n{k+l-1}E_{k+l-1}(x)E_{n-k}(y)
\\=&\sum_{k=0}^n\bi nk\bi{k+n}{k+l}\l((-1)^{n-k}B_{k+l}(x)-B_{k+l}(y)\r)E_{n-k}(x-y)
\endaligned\tag1.6$$
$($where $\da_{l,m}$ takes $1$ or $0$ according to whether $l=m$ or not$)$, and
$$\aligned&\f{n-l}n\sum_{k=0}^{n-l}\bi nk\bi n{k+l}B_{k+l}(x)B_{n-k}(y)
\\=&\sum_{k=0}^n\bi nk\bi{k+n-1}{k+l}\l((-1)^{n-k}B_{k+l}(x)+B_{k+l}(y)\r)B_{n-k}(x-y).
\endaligned\tag1.7$$
In particular,
$$\aligned&\f{(n+1)(n+1-l)}8\sum_{k=\da_{l,0}}^n\bi nk\bi n{k+l-1}E_{k+l-1}(x)E_{n-k}(x)
\\&=\sum_{k=0}^{n-1}\bi{n+1}k\bi{k+n}{k+l}B_{k+l}(x)\l(2^{n-k+1}-1\r)B_{n-k+1}
\endaligned\tag1.8$$
and
$$\aligned&\sum_{k=0}^{n-l}\bi nk\bi n{k+l}B_{k+l}(x)B_{n-k}(x)
\\=&\f{2n}{n-l}\sum^n\Sb k=0\\k\not=n-1\endSb\bi nk\bi {k+n-1}{k+l}B_{k+l}(x)B_{n-k}.
\endaligned\tag1.9$$
\endproclaim

\Proof. As $(l-n-1)+n+n=(n+l)-1$ and $(1-x)+y+(x-y)=1$, by Theorem 1.1(i) we have
$$\align&\sum_{k=0}^{n+l}(-1)^k\bi{l-n-1}k\bi n{n+l-k}B_k(1-x)E_{n+l-k}(x-y)
\\&-(-1)^{n+l}\sum_{k=0}^{n+l}(-1)^k\bi{l-n-1}k\bi n{n+l-k}B_k(y)E_{n+l-k}(x-y)
\\=&\f{l-n-1}2\sum_{k=0}^{n-\da_{l,0}}(-1)^k\bi nk\bi n{n+l-1-k}E_k(y)E_{n+l-1-k}(1-x)
\\=&\f{l-n-1}2\sum_{k=\da_{l,0}}^n(-1)^{n-k}\bi nk\bi n{k+l-1}E_{n-k}(y)E_{k+l-1}(1-x)
\endalign$$
which can be reduced to (1.6). (1.8) follows from (1.6) in the case $y=x$ since
$((-1)^m-1)E_m(0)=4(2^{m+1}-1)B_{m+1}/(m+1)$ for $m=1,2,3,\ldots$.
(It is known that $(m+1)E_m(x)=2(B_{m+1}(x)-2^{m+1}B_{m+1}(x/2))$ (cf. [AS] and [S1]).)

In light of Theorem 1.1(ii),
$$(l-n)\mm nn{1-x}y_{n+l}+n\mm n{l-n}y{x-y}_{n+l}+n\mm{l-n}n{x-y}{1-x}_{n+l}=0.$$
This is equivalent to (1.7). In the case $y=x$, (1.7) gives (1.9)
because $((-1)^{m}+1)B_m=2B_m$ for $m=0,2,3,\ldots$. \qed

\Remark\ 1.3. Putting $l=0$ and $x=1/2$ in (1.8)
and noting that $B_k(1/2)=(2^{1-k}-1)B_k$ (see, e.g., [AS] and [S1]), we then
get the following identity:
$$\aligned&\f{(n+1)^2}8\sum_{k=0}^{n-1}\bi nk\bi n{k+1}E_kE_{n-1-k}
\\=&-\sum_{k=0}^{n-1}\bi{n+1}k\bi{n+k}n2^{n-k}(2^{k-1}-1)(2^{n-k+1}-1)B_kB_{n-k+1}
\endaligned$$
for any $n\in\Z^+$. Similarly, (1.9) in the case $l=x=0$ yields the following new identity:
$$\sum_{k=2}^{n-2}\bi nk^2B_kB_{n-k}
-2\sum_{k=2}^{n-2}\bi nk\bi{n+k-1}kB_kB_{n-k}=2\bi{2n-1}{n-1}B_n$$
for every $n=4,5,\ldots$.
\medskip

The following theorem can be deduced from Theorem 1.1.

\proclaim{Theorem 1.2} Let $l,m,n\in\Z^+$, $l\ls \min\{m,n\}$ and $x+y+z=1$. Then
$$\aligned&(-1)^m\sum_{k=0}^m\bi mk\bi{n+k}{l-1}B_{n-l+k+1}(x)E_{m-k}(z)
\\&+(-1)^{n-l}\sum_{k=0}^n\bi nk\bi{m+k}{l-1}B_{m-l+k+1}(y)E_{n-k}(z)
\\=&-\f l2\sum_{k=0}^l(-1)^k\bi mk\bi n{l-k}E_{n-l+k}(x)E_{m-k}(y)
\endaligned\tag1.10$$
and
$$\aligned&\sum_{k=0}^{l}(-1)^k\bi mk\bi n{l-k}B_{m-k}(x)E_{n-l+k}(z)
\\&-(-1)^m\sum_{k=0}^m\bi mk\bi{n+k}{l}B_{m-k}(y)E_{n-l+k}(z)
\\=&(-1)^{n-l-1}\f m2\sum^n_{k=\da_{l,m}}\bi nk\bi{m+k-1}{l}E_{n-k}(y)E_{m-l-1+k}(x).
\endaligned\tag1.11$$
We also have
$$\aligned&\f{(-1)^m}m\sum_{k=0}^m\bi mk\bi{n+k-1}{l-1}B_{n-l+k}(x)B_{m-k}(z)
\\&+(-1)^l\f{(-1)^n}n\sum_{k=0}^n\bi nk\bi{m+k-1}{l-1}B_{m-l+k}(y)B_{n-k}(z)
\\&=\f l{mn}\sum_{k=0}^l(-1)^k\bi mk\bi n{l-k}B_{n-l+k}(x)B_{m-k}(y).
\endaligned\tag1.12$$
\endproclaim

\proclaim{Corollary 1.3 {\rm (Woodcock [W])}} Let $m,n\in\Z^+$. Then
$$\f1m\sum_{k=1}^m\bi mk(-1)^kB_{m-k}B_{n-1+k}
=\f1n\sum_{k=1}^n\bi nk(-1)^kB_{n-k}B_{m-1+k}.$$
\endproclaim
\Proof. Simply take $x=y=0$ and $l=z=1$ in (1.12). \qed

From Theorem 1.1 we can also deduce the following result.

\proclaim{Theorem 1.3} Let $n\in\Z^+$, and let $t,x,y,z$ be parameters with $x+y+z=1$.
Then we have
$$\aligned& \f{(-1)^n}2\sum_{k=0}^{n-1}\bi t{k}E_k(x)E_{n-1-k}(y)
\\=&\f 1{n-t}\sum_{k=0}^n\bi {n-t}kB_k(x)E_{n-k}(z)
+\bi tn\sum_{k=0}^n\bi nk\f {E_{k}(z)}{t-k}B_{n-k}(y)
\endaligned\tag1.13$$
and
$$\aligned&\f n2\bi tn\sum_{k=0}^{n-1}\bi {n-1}k\f {E_{k}(x)}{t-k}E_{n-1-k}(y)
-(-1)^nE_n(z)\bi tn\sum_{k=0}^{n-1}\f1{t-k}
\\=&(-1)^n\sum_{k=1}^n\bi t{n-k}\f {B_k(y)}{k}E_{n-k}(z)
-\sum_{k=1}^n\bi {n-1-t}{n-k}\f {B_k(x)}kE_{n-k}(z).
\endaligned\tag1.14$$
Also,
$$\aligned&\f {(-1)^{n-1}}n\bi {t-1}{n-1}\sum_{k=0}^n\bi nk\f {B_k(x)}{t-k}B_{n-k}(y)
-\f {B_n(z)}n\bi {t-1}{n-1}\sum_{k=1}^{n-1}\f1{t-k}
\\&=\f 1t\sum_{k=1}^{n}\bi t{n-k}\f{B_{k}(y)}kB_{n-k}(z)
+\f {(-1)^n}{n-t}\sum_{k=1}^n\bi {n-t}{n-k}\f{B_k(x)}kB_{n-k}(z).
\endaligned\tag1.15$$
\endproclaim

\proclaim{Corollary 1.4} Let $n\in\Z^+$ and $x+y+z=1$. Then
$$\aligned&\sum_{k=0}^n\bi{n+1}k\l((-1)^nB_k(x)-B_k(y)\r)E_{n-k}(z)
\\&=\f{n+1}2\sum_{l=0}^{n-1}(-1)^lE_l(x)E_{n-1-l}(y),\endaligned\tag1.16$$
$$\aligned &\sum_{k=1}^n\bi nk\f{B_k(x)}kE_{n-k}(z)-\sum_{k=1}^n(-1)^k\f{B_k(y)}kE_{n-k}(z)
\\&=\f{(-1)^n}2\sum_{l=0}^{n-1}\bi nlE_l(y)E_{n-1-l}(x)-H_nE_n(z)
\endaligned\tag1.17$$
and
$$\aligned&(-1)^n\sum_{k=0}^n\bi{n+1}kB_{n-k}(x)B_k(y)
+\sum_{k=0}^{n-1}\bi {n+1}k\f{B_{n-k}(x)}{n-k}B_k(z)
\\&=(n+1)\sum_{k=1}^n(-1)^k\f{B_k(y)}kB_{n-k}(z)+(1-H_n)(n+1)B_n(z).
\endaligned\tag1.18$$
\endproclaim
\Proof. Taking $t=-1$ in Theorem 1.3 we immediately get (1.16)-(1.18). \qed

\proclaim{Corollary 1.5} Let $n\in\Z^+$ and $x+y+z=1$. Then
$$\aligned&\f12\sum_{k=1}^{n-1}(-1)^{k-1}\f{E_k(x)}kE_{n-1-k}(y)+\f{H_{n-1}E_{n-1}(y)}2
\\=&\f1n\sum_{k=1}^n\bi nk\f{E_k(z)}kB_{n-k}(y)
+\f{(-1)^n}n\sum_{k=1}^n\bi nkH_kE_k(z)B_{n-k}(x)
\endaligned\tag1.19$$
and
$$\aligned&\f{(-1)^{n-1}}2\sum_{k=1}^{n-1}\bi{n-1}k\f{E_k(x)}k
E_{n-1-k}(y)+H_{n-1}\f{E_n(z)+(-1)^nB_n(y)}n
\\&=\sum_{k=1}^{n-1}(-1)^k\f{B_k(y)}k\cdot\f{E_{n-k}(z)}{n-k}
+\sum_{k=1}^n\bi{n-1}{k-1}H_{k-1}\f{B_k(x)}kE_{n-k}(z).
\endaligned\tag1.20$$
We also have
$$\aligned&\sum_{k=1}^n\bi {n-1}{k-1}\f{B_k(x)}{k^2}\l(B_{n-k}(y)+(-1)^nB_{n-k}(z)\r)
\\=&\sum_{k=1}^{n-1}(-1)^{n-k}\f{B_k(y)}k\cdot\f{B_{n-k}(z)}{n-k}-H_{n-1}
\f{B_n(y)+(-1)^nB_n(z)}n
\endaligned\tag1.21$$
\endproclaim

\Remark\ 1.4. In the case $x=y=0$ and $z=1$, (1.21) yields Miki's identity.

\medskip

The next section is devoted to proofs of Theorems 1.1 and 1.2.
Theorem 1.3 and Corollary 1.5 will be proved in Section 3.

\heading{2. Proofs of Theorems 1.1--1.2}\endheading

\proclaim{Lemma 2.1} Let $P(x),Q(x)\in \C[x]$
where $\C$ is the field of complex numbers.

{\rm (i)} We have
$$\Delta(P(x)Q(x))=P(x+1)\Delta(Q(x))+\Delta(P(x))Q(x)\tag2.1$$
and
$$\Delta^*(P(x)Q(x))=P(x+1)\Delta^*(Q(x))-\Delta(P(x))Q(x).\tag2.2$$

{\rm (ii) If $\Delta(P(x))=\Delta(Q(x))$, then $P'(x)=Q'(x)$.
 If $\Delta^*(P(x))=\Delta^*(Q(x))$, then $P(x)=Q(x)$.
\endproclaim
\Proof. The first part can be verified easily. Part
(ii) is Lemma 3.1 of [PS]. \qed

\smallskip

The following lemma has the same flavor with Theorem 1.1 of Sun [S2].

\proclaim{Lemma 2.2} Let $\{a_l\}_{l=0}^{\infty}$ be a sequence of
complex numbers, and $\{a_l^*\}_{l=0}^{\infty}$ be its dual sequence.
Set
$$A_k(t)=\sum_{l=0}^k\bi kl(-1)^la_lt^{k-l}\ \t{and}\ A_k^*(t)
=\sum_{l=0}^k\bi kl(-1)^la_l^*t^{k-l}\tag2.3$$
for $k=0,1,2,\ldots$. Let $n\in\Z^+$,
$r+s+t=n-1$ and $x+y+z=1$. Then
$$\sum_{k=0}^n(-1)^k\bi rkx^{n-k}\(\bi {s}{n-k}A_k(y)-(-1)^n\bi t{n-k}A_k^*(z)\)=0.\tag2.4$$
\endproclaim
\Proof. By Remark 1.1 of Sun [S2],
$$(-1)^kA_k^*(z)=A_k(x+y)=\sum_{l=0}^k\bi klx^{k-l}A_l(y).$$
Therefore
$$\align&\sum_{k=0}^n(-1)^k\bi rk\bi t{n-k}x^{n-k}A_k^*(z)
\\=&\sum_{k=0}^n\bi rk\bi t{n-k}x^{n-k}\sum_{l=0}^k\bi klx^{k-l}A_l(y)
\\=&\sum_{l=0}^nx^{n-l}A_l(y)\sum_{k=l}^n\bi rl\bi{r-l}{k-l}\bi t{n-k}
\\=&\sum_{l=0}^n\bi rlx^{n-l}A_l(y)c_l
\endalign$$
where
$$\align c_l=&\sum_{k=l}^n\bi{r-l}{k-l}\bi t{n-k}=\bi{r+t-l}{n-l}
\ (\t{by Vandermonde's identity})
\\=&(-1)^{n-l}\bi {l-r-t+n-l-1}{n-l}=(-1)^{n-l}\bi s{n-l}.
\endalign$$
Thus (2.4) follows. \qed

\Remark\ 2.1. If we let $a_l=(-1)^lB_l$ for $l=0,1,2,\ldots$, then
$A_k(t)=A_k^*(t)=B_k(t)$. Also, $A_k(t)=A_k^*(t)=E_k(t)$
if $a_l=(-1)^lE_l(0)$ for $l=0,1,2,\ldots$.
\medskip

\noindent
{\tt Proof of Theorem 1.1}. We fix $y$ and view $z=1-x-y$ as a function in $x$.

(i) Set $$P(x)=\sum_{k=0}^n(-1)^k\bi rk\bi s{n-k}B_k(x)E_{n-k}(z).$$
Then, by Lemma 2.1, $\Delta^*(P(x))$ coincides with
$$\align&\sum_{k=0}^n(-1)^k\bi rk\bi s{n-k}\Delta^*(B_k(x)E_{n-k}(z))
\\=&\sum_{k=0}^n(-1)^k\bi rk\bi s{n-k}\l(B_k(x+1)2(z-1)^{n-k}-kx^{k-1}E_{n-k}(z)\r)
\\=&2\sum_{k=0}^n(-1)^k\bi rk\bi s{n-k}(z-1)^{n-k}B_k(x+1)+r\Sigma
\endalign$$
where
$$\align\Sigma=&\sum_{k=1}^n(-1)^{k-1}\bi {r-1}{k-1}\bi s{n-k}x^{k-1}E_{n-k}(z).
\\=&(-1)^{n-1}\sum_{l=0}^{n-1}(-1)^l\bi {r-1}{n-1-l}\bi slx^{n-1-l}E_l(z).
\endalign$$
Applying Lemma 2.2 and Remark 2.1 we obtain that
$$\align \Delta^*(P(x))=&2(-1)^n\sum_{k=0}^n(-1)^k\bi rk\bi t{n-k}(z-1)^{n-k}B_k(y)
\\&+r\sum_{l=0}^{n-1}(-1)^l\bi sl\bi t{n-1-l}x^{n-1-l}E_l(y).
\endalign$$
It follows that $\Delta^*(P(x))=\Delta^*(Q(x))$ where
$$\align Q(x)=&(-1)^n\sum_{k=0}^n(-1)^k\bi rk\bi t{n-k}B_k(y)E_{n-k}(z)
\\&+\f r2\sum_{l=0}^{n-1}(-1)^l\bi sl\bi t{n-1-l} E_l(y)E_{n-1-l}(x).
\endalign$$
Thus $P(x)=Q(x)$ by Lemma 2.1. This is equivalent to the desired (1.1).
\qed

(ii) Set $$P_n(x)=\M rszx=\sum_{k=0}^n(-1)^k\bi rk\bi s{n-k}B_k(x)B_{n-k}(z).$$
By Lemma 2.1,
$$\align \Delta(B_k(x)B_{n-k}(z))=&\Delta(B_k(x))B_{n-k}(z)+B_k(x+1)\Delta(B_{n-k}(z))
\\=&kx^{k-1}B_{n-k}(z)-(n-k)B_k(x+1)(z-1)^{n-k-1}
\endalign$$
for every $k=0,1,\ldots,n$. Thus
$$\Delta(P_n(x))=r R(x)-s\sum_{k=0}^{n-1}(-1)^k\bi rk\bi {s-1}{n-k-1}B_k(x+1)(z-1)^{n-k-1}$$
where
$$\align R(x)=&\sum_{k=1}^n(-1)^k\bi {r-1}{k-1}\bi s{n-k}x^{k-1}B_{n-k}(z)
\\=&(-1)^n\sum_{l=0}^{n-1}(-1)^l\bi sl\bi {r-1}{n-1-l}x^{n-1-l}B_l(z).
\endalign$$
Applying Lemma 2.2 and Remark 2.1 we obtain that
$$\align \Delta(P_n(x))=&-r\sum_{l=0}^{n-1}(-1)^l\bi sl\bi {t-1}{n-1-l}x^{n-1-l}B_l(y)
\\&-s(-1)^{n-1}\sum_{l=0}^{n-1}(-1)^l\bi rl\bi {t-1}{n-1-l}(z-1)^{n-1-l}B_l(y)
\endalign$$
It follows that $\Delta(P_n(x))=\Delta(Q_n(x))$ where
$$\align Q_n(x)=&-\f rt\sum_{l=0}^{n-1}(-1)^l\bi sl\bi {t}{n-l}B_{n-l}(x)B_l(y)
\\&-(-1)^n\f st\sum_{l=0}^{n-1}(-1)^l\bi rl\bi t{n-l}B_{n-l}(z)B_l(y).
\\=&-\f rt\sum_{l=0}^{n-1}(-1)^l\bi sl\bi {t}{n-l}B_{n-l}(x)B_l(y)
\\&-\f st\sum_{k=1}^n(-1)^k\bi t{k}\bi r{n-k}B_{k}(z)B_{n-k}(y).
\endalign$$
Thus $P_n'(x)=Q_n'(x)$ by Lemma 2.1.

Observe that $P'_n(x)$ coincides with
$$\align &\sum_{k=1}^n(-1)^k\bi rk\bi s{n-k}kB_{k-1}(x)B_{n-k}(z)
\\&-\sum_{k=0}^{n-1}(-1)^k\bi rk\bi s{n-k}(n-k)B_k(x)B_{n-k-1}(z)
\\=&\sum_{k=0}^{n-1}(-1)^{k+1}\bi r{k+1}\bi s{n-1-k}(k+1)B_k(x)B_{n-1-k}(z)
\\&-\sum_{k=0}^{n-1}(-1)^k\bi rk\bi s{n-k}(n-k)B_k(x)B_{n-1-k}(z)
\\=&\sum_{k=0}^{n-1}(-1)^{k-1}\bi rk\bi s{n-1-k}(r-k+(s-n+k+1))B_k(x)B_{n-1-k}(z)
\\=&(t-1)\m rszx
\endalign$$
and
$$\align Q_n'(x)=&-r\sum_{l=0}^{n-1}(-1)^l\bi sl\bi{t-1}{n-l-1}B_{n-l-1}(x)B_l(y)
\\&+s\sum_{k=1}^n(-1)^k\bi{t-1}{k-1}\bi r{n-k}B_{k-1}(z)B_{n-k}(y)
\\=&-r\m s{t-1}xy-s\m {t-1}ryz.
\endalign$$
Thus the equality $P_n'(x)=Q_n'(x)$ gives that
$$r\m s{t'}xy+s\m {t'}ryz+t'\m rszx=0$$
where $t'=t-1=n-1-(r+s)$. Replacing $n-1$ by $n$ we then obtain the required identity
(1.2). This concludes the proof. \qed

\medskip
\noindent{\tt Proof of Theorem 1.2}. Clearly $\bar n=m+n-l\in\Z^+$. By Theorem 1.1(i),
$$\align&\sum_{k=0}^{\bar n+1}(-1)^k\bi{-l}k\bi m{\bar n+1-k}B_k(x)E_{\bar n+1-k}(z)
\\&-(-1)^{\bar n+1}\sum_{k=0}^{\bar n+1}(-1)^k\bi{-l}k\bi n{\bar n+1-k}B_k(y)E_{\bar n+1-k}(z)
\\=&\f{-l}2\sum_{k=0}^{\bar n}(-1)^k\bi mk\bi n{\bar n-k}E_k(y)E_{\bar n-k}(x).
\endalign$$
That is,
$$\align&\sum_{k=0}^m(-1)^{\bar n+1-k}\bi{-l}{\bar n+1-k}\bi mkB_{\bar n+1-k}(x)E_k(z)
\\&-\sum_{k=0}^n(-1)^k\bi{-l}{\bar n+1-k}\bi nkB_{\bar n+1-k}(y)E_k(z)
\\=&\f{-l}2\sum_{k=0}^m(-1)^{m-k}\bi mk\bi n{n-l+k}E_{m-k}(y)E_{n-l+k}(x)
\\=&(-1)^{m-1}\f l2\sum_{k=0}^l(-1)^k\bi mk\bi n{l-k}E_{m-k}(y)E_{n-l+k}(x).
\endalign$$
Therefore (1.10) follows.
By Theorem 1.1(i) we also have
$$\align&\sum_{k=0}^{\bar n}(-1)^k\bi mk\bi n{\bar n-k}B_k(x)E_{\bar n-k}(z)
\\&-(-1)^{\bar n}\sum_{k=0}^{\bar n}(-1)^k\bi mk\bi{-l-1}{\bar n-k}B_k(y)E_{\bar n-k}(z)
\\=&\f m2\sum_{k=0}^{\bar n-1}(-1)^k\bi nk\bi{-l-1}{\bar n-1-k}E_k(y)E_{\bar n-1-k}(x)
\\=&\f m2\sum_{k=0}^{n-\da_{l,m}}(-1)^k\bi nk\bi{-l-1}{m+n-l-1-k}E_k(y)E_{m+n-l-1-k}(x)
\\=&\f m2\sum_{k=\da_{l,m}}^n(-1)^{n-k}\bi nk\bi{-l-1}{m-l-1+k}E_{n-k}(y)E_{m-l-1+k}(x),
\endalign$$
which gives (1.11) after few trivial steps.

In light of Theorem 1.1(ii),
$$l\mm mnxy_{\bar n}=m\mm n{-l}yz_{\bar n}+n\mm{-l}mzx_{\bar n}.$$
That is,
$$\align&l\sum_{k=0}^m(-1)^{m-k}\bi mk\bi n{n-l+k}B_{n-l+k}(x)B_{m-k}(y)
\\=&m\sum_{k=0}^n(-1)^{n-k}\bi nk\bi{-l}{m-l+k}B_{m-l+k}(y)B_{n-k}(z)
\\&+n\sum_{k=0}^m(-1)^{\bar n-k}\bi{-l}{\bar n-k}\bi mkB_k(z)B_{\bar n-k}(x).
\endalign$$
This is equivalent to (1.12). We are done. \qed

\heading{3. Proofs of Theorem 1.3 and Corollary 1.5}\endheading

\proclaim{Lemma 3.1} Let $n$ be a nonnegative integer and $s$ be a parameter. Then
$$\lim_{t\to0}\f1t\(\bi{s+t}n-\bi sn\)=\bi sn\sum_{0\ls l<n}\f1{s-l}.\tag3.1$$
In particular,
$$\lim_{t\to 0}\f 1t\(\bi{t-1}n-(-1)^n\)=(-1)^{n-1}H_n.\tag3.2$$
\endproclaim
\Proof. Observe that
$$\bi{s+t}n=\bi sn\prod_{0\ls l<n}\f{s+t-l}{s-l}=\bi sn\prod_{0\ls l<n}\l(1+\f t{s-l}\r).$$
So (3.1) follows. In the case $s=-1$, (3.1) turns out to be (3.2). \qed

\medskip
\noindent
{\tt Proof of Theorem 1.3}. (1.1) in the case $s=-1$ yields that
$$\aligned&(-1)^n\sum_{k=0}^n\bi {n-t}kB_k(x)E_{n-k}(z)
\\&-(-1)^n\sum_{k=0}^n(-1)^k\bi {n-t}k\bi t{n-k}B_k(y)E_{n-k}(z)
\\=&\f {n-t}2\sum_{l=0}^{n-1}\bi t{n-1-l}E_{n-1-l}(x)E_l(y).
\endaligned$$
For each $k=0,1,\ldots,n$ we clearly have
$$\aligned\bi {n-t}k\bi t{n-k}=&\bi nk\bi tn\f{(n-t)(n-t-1)
\cdots(n-t-k+1)}{(t-n+k)\cdots(t-n+1)}
\\=&(-1)^k\bi nk\bi tn\f {t-n}{t-n+k}.
\endaligned$$
Therefore
$$\aligned&\f{(-1)^n}2\sum_{k=0}^{n-1}\bi t{k}E_k(x)E_{n-1-k}(y)
-\f 1{n-t}\sum_{k=0}^n\bi {n-t}kB_k(x)E_{n-k}(z)
\\=&\bi tn\sum_{k=0}^n\bi nk\f {B_k(y)}{t+k-n}E_{n-k}(z)
=\bi tn\sum_{l=0}^n\bi nl\f {E_{l}(z)}{t-l}B_{n-l}(y).
\endaligned$$
This proves (1.13).

Now we come to prove (1.14) and view $s=n-1-r-t$ as a function in $r$.
In light of (1.1),
$$\align&\f 12\sum_{l=0}^{n-1}(-1)^l\bi sl\bi t{n-1-l}E_l(y)E_{n-1-l}(x)
\\=&\f1r\sum_{k=0}^n(-1)^k\bi rkE_{n-k}(z)\(\bi s{n-k}B_k(x)-(-1)^n\bi t{n-k}B_k(y)\)
\\=&\sum_{k=1}^n\f{(-1)^k}k\bi{r-1}{k-1}E_{n-k}(z)\(\bi s{n-k}B_k(x)-(-1)^n\bi t{n-k}B_k(y)\)
\\&+(-1)^nE_n(z)\f{(-1)^n\bi sn-\bi tn}r.
\endalign$$
By Lemma 3.1,
$$\lim_{r\to 0}\f 1r\((-1)^n\bi sn-\bi t{n}\)=\lim_{r\to 0}\f 1r\(\bi{r+t}n-\bi t{n}\)
=\bi tn\sum_{l=0}^{n-1}\f1{t-l}.$$ As in the proof of (1.13), we
also have
$$\align(-1)^l\bi {n-1-t}l\bi t{n-1-l}=&\bi{n-1}l\bi t{n-1}\f{t-(n-1)}{t-(n-1)+l}
\\=&\f n{t+l-(n-1)}\bi tn\bi{n-1}l
\endalign$$
for every $l=0,1,\ldots,n-1$.
Thus, by letting $r\to 0$ we get from the above that
$$\align&\f n2\bi tn\sum_{l=0}^{n-1}\bi{n-1}l\f{E_l(y)E_{n-1-l}(x)}{t+l-n+1}
-(-1)^nE_n(z)\bi tn\sum_{l=0}^{n-1}\f1{t-l}
\\&=-\sum_{k=1}^nE_{n-k}(z)\(\bi {n-1-t}{n-k}\f{B_k(x)}k-(-1)^n\bi t{n-k}\f{B_k(y)}k\),
\endalign$$
which is equivalent to (1.14).

 Now we turn to prove (1.15). Let us view $s=n-r-t$ as a function in $r$.
 Then
 $$\align\lim_{r\to 0}\M stxy=&\M{n-t}txy=\sum_{k=0}^n\bi nk\bi tn\f{t-n}{t-n+k}B_{n-k}(x)B_k(y)
 \\=&(t-n)\bi tn\sum_{l=0}^n\bi nl\f{B_l(x)}{t-l}B_{n-l}(y).
 \endalign$$
 On the other hand,
 $$\align &\lim_{r\to 0}\f 1r\(s\M tryz+t\M rszx\)
 \\=&(n-t)(-1)^{n-1}\sum_{k=0}^{n-1}\bi tk\f{B_{n-k}(y)}{n-k}B_k(z)
 \\&-t\sum_{k=1}^n\bi{n-t}{n-k}\f{B_k(x)}kB_{n-k}(z)+(-1)^nB_n(z)R
 \endalign$$
 where
 $$\align R=&\lim_{r\to 0}\f1r\((n-t-r)\bi tn+(-1)^nt\bi{n-t-r}n\)
 \\=&\lim_{r\to 0}\f1r\(t\bi{r+t-1}n-(t-n)\bi tn\)-\bi tn
 \\=&\lim_{r\to 0}\f tr\(\bi {r+t-1}n-\bi{t-1}n\)-\bi tn
 \\=&t\bi{t-1}n\sum_{l=0}^{n-1}\f1{t-1-l}-\bi tn=t\bi{t-1}n\sum_{k=1}^{n-1}\f1{t-k}.
 \endalign$$
 Applying (1.2) we then get (1.15) from the above.

 The proof of Theorem 1.3 is now complete. \qed

\medskip
\noindent{\tt Proof of Corollary 1.5}.
We can easily get (1.21) by calculating the limitation of
the left hand side of (1.15) minus the right hand side of (1.15)
as $t$ tends to $0$. Thus it remains to show (1.19) and (1.20).

(1.13) can be rewritten in the form
$$\aligned&\f{(-1)^n}2t\sum_{k=1}^{n-1}\bi{t-1}{k-1}\f{E_k(x)}kE_{n-1-k}(y)
+\f{(-1)^n}2E_{n-1}(y)
\\=&\sum_{k=0}^n\(\f{\bi{n-t}k}{n-t}-\f{\bi nk}{n}\)B_k(x)E_{n-k}(z)
+\f 1n\sum_{k=0}^n\bi nkB_k(x)E_{n-k}(z)
\\&+\f tn\bi{t-1}{n-1}\(\f{B_n(y)}t+\sum_{k=1}^n\bi nk\f{E_k(z)}{t-k}B_{n-k}(y)\).
\endaligned$$
Letting $t\to 0$ we get that
$$\f 1n\sum_{k=0}^n\bi nkB_k(x)E_{n-k}(z)+(-1)^{n-1}\f{B_n(y)}n
=\f{(-1)^n}2E_{n-1}(y).\tag3.3$$
Thus
$$\aligned&\f{(-1)^n}2\sum_{k=1}^{n-1}\bi{t-1}{k-1}\f{E_k(x)}kE_{n-1-k}(y)
\\=&\sum_{k=0}^n\f 1t\(\f{\bi{n-t}k}{n-t}-\f{\bi nk}{n}\)B_k(x)E_{n-k}(z)
+\f{B_n(y)}{nt}\(\bi{t-1}{n-1}-(-1)^{n-1}\)
\\&+\f 1n\bi{t-1}{n-1}\sum_{k=1}^n\bi nk\f{E_k(z)}{t-k}B_{n-k}(y).
\endaligned$$
Letting $t\to 0$ we then have
$$\aligned&\f{(-1)^n}2\sum_{k=1}^{n-1}\f{(-1)^{k-1}}kE_k(x)E_{n-1-k}(y)
+\f{(-1)^{n-1}}n\sum_{k=1}^n\bi nk\f{E_k(z)}kB_{n-k}(y)
\\=&\sum_{k=0}^n\lim_{t\to 0}\f{n\bi{n-t}k-(n-t)\bi nk}{tn(n-t)}B_k(x)E_{n-k}(z)
+\f{B_n(y)}n(-1)^nH_{n-1}.
\endaligned$$
Observe that
$$\aligned&\lim_{t\to 0}\f{n\bi{n-t}k-(n-t)\bi nk}{t(n-t)}
=\lim_{t\to 0}\(\f{\bi nk}{n-t}-\f n{n-t}\cdot\f{\bi{n-t}k-\bi nk}{-t}\)
\\=&\f 1n\bi nk-\bi nk\sum_{l=0}^{k-1}\f 1{n-l}
=-\bi nk\sum_{0<l<k}\f 1{n-l}=\bi nk(H_{n-k}-H_{n-1}).
\endaligned$$
Therefore
$$\aligned&\f{(-1)^{n-1}}2\sum_{k=1}^{n-1}\f{(-1)^k}kE_k(x)E_{n-1-k}(y)+
\f{(-1)^{n-1}}n\sum_{k=1}^n\bi nk\f{E_k(z)}kB_{n-k}(y)
\\=&(-1)^nH_{n-1}\f{B_n(y)}n+\f 1n\sum_{k=0}^n\bi nk(H_{n-k}-H_{n-1})B_k(x)E_{n-k}(z)
\\=&(-1)^nH_{n-1}\f{B_n(y)}n-\f{H_{n-1}}n\sum_{k=0}^n\bi nkB_k(x)E_{n-k}(z)
\\&+\f 1n\sum_{l=0}^n\bi nlH_lE_l(z)B_{n-l}(x)
\\=&-H_{n-1}\f{(-1)^n}2E_{n-1}(y)+\f 1n\sum_{k=0}^n\bi nkH_kE_k(z)B_{n-k}(x).
\endaligned$$
This proves (1.19).

We can reformulate (1.14) as follows:
$$\aligned&\f t2\bi{t-1}{n-1}\sum_{k=1}^{n-1}\bi{n-1}k\f{E_k(x)}{t-k}E_{n-1-k}(y)
+\f 12\bi{t-1}{n-1}E_{n-1}(y)
\\&-(-1)^nE_n(z)\f tn\bi{t-1}{n-1}\sum_{k=1}^{n-1}\f 1{t-k}-(-1)^n\f{E_n(z)}n\bi{t-1}{n-1}
\\=&(-1)^nt\sum_{k=1}^{n-1}\bi{t-1}{n-k-1}\f{B_k(y)}k\cdot\f{E_{n-k}(z)}{n-k}+(-1)^n\f{B_n(y)}n
\\&-\sum_{k=1}^n\(\bi{n-1-t}{n-k}-\bi{n-1}{n-k}\)\f{B_k(x)}kE_{n-k}(z)
\\&-\sum_{k=1}^n\bi{n-1}{n-k}\f{B_k(x)}kE_{n-k}(z).
\endaligned$$
In view of (3.3),
$$\align\sum_{k=1}^n\bi{n-1}{n-k}\f{B_k(x)}kE_{n-k}(z)
=&\sum_{k=1}^n\bi{n-1}{k-1}\f{B_k(x)}kE_{n-k}(z)
\\=&(-1)^n\l(\f{B_n(y)}n+\f{E_{n-1}(y)}2\r)-\f{E_n(z)}n.
\endalign$$
Therefore
$$\aligned&\f 12\bi{t-1}{n-1}\sum_{k=1}^{n-1}\bi{n-1}k\f{E_k(x)}{t-k}E_{n-1-k}(y)
\\&-(-1)^n\f{E_n(z)}n\bi{t-1}{n-1}\sum_{k=1}^{n-1}\f 1{t-k}
\\=&(-1)^n\sum_{k=1}^{n-1}\bi{t-1}{n-k-1}\f{B_k(y)}k\cdot\f{E_{n-k}(z)}{n-k}
\\&+\sum_{k=1}^{n-1}\f{\bi{n-1-t}{n-k}-\bi{n-1}{n-k}}{-t}\cdot\f{B_k(x)}kE_{n-k}(z)
\\&-\l(\f{E_{n-1}(y)}2+(-1)^{n-1}\f{E_n(z)}n\r)\f{\bi{t-1}{n-1}-(-1)^{n-1}}t.
\endaligned$$
Letting $t\to 0$ we obtain that
$$\align&\f{(-1)^n}2\sum_{k=1}^{n-1}\bi{n-1}k\f{E_k(x)}kE_{n-1-k}(y)-\f{E_n(z)}nH_{n-1}
\\=&\sum_{k=1}^{n-1}(-1)^{k-1}\f{B_k(y)}k\cdot\f{E_{n-k}(z)}{n-k}
\\&+\sum_{k=1}^{n-1}\bi{n-1}{n-k}\(\sum_{l=0}^{n-k-1}\f 1{n-1-l}\)\f{B_k(x)}kE_{n-k}(z)
\\&+H_{n-1}\(\f{(-1)^{n-1}}2E_{n-1}(y)+\f{E_n(z)}n\).
\endalign$$
It follows that
$$\align &\f{(-1)^n}2\sum_{k=1}^{n-1}\bi{n-1}k\f{E_k(x)}kE_{n-1-k}(y)
+\sum_{k=1}^{n-1}(-1)^{k}\f{B_k(y)}k\cdot\f{E_{n-k}(z)}{n-k}
\\=&\sum_{k=1}^n\bi{n-1}{n-k}(H_{n-1}-H_{k-1})\f{B_k(x)}kE_{n-k}(z)
\\&+H_{n-1}\(\f{(-1)^{n-1}}2E_{n-1}(y)+2\f{E_n(z)}n\)
\\=&-\sum_{k=1}^n\bi{n-1}{k-1}H_{k-1}\f{B_k(x)}kE_{n-k}(z)+H_{n-1}R
\endalign$$
where
$$\align R=&\sum_{k=1}^n\bi{n-1}{k-1}\f{B_k(x)}kE_{n-k}(z)+\f{(-1)^{n-1}}2E_{n-1}(y)+\f 2nE_n(z)
\\=&\f1n\l(E_n(z)+(-1)^nB_n(y)\r)\ \quad\ (\t{by (3.3)}).
\endalign$$
This proves (1.20). We are done. \qed

\Ack. The authors thank the referee for his/her helpful comments.

\bigskip

\leftline{Department of Mathematics} \leftline{Nanjing University}
\leftline{Nanjing 210093} \leftline{People's Republic of China}
\leftline {E-mail: (Zhi-Wei Sun) {\tt zwsun\@nju.edu.cn}}
\leftline {\quad\qquad\ \ (Hao Pan) {\tt haopan79\@yahoo.com.cn}}

\enddocument